\documentclass[a4paper,12pt,leqno]{amsart}
\usepackage{amsmath,amssymb,mathrsfs}
\usepackage[matrix,arrow]{xy} 
\usepackage{xcolor,hyperref}
\hypersetup{colorlinks=true,citecolor=black,linkcolor=black}
\newcommand{\h}{\hbox}
\newcommand{\q}{\quad}
\newcommand{\nin}{\noindent}

\newcommand{\sk}{\par\smallskip}

\newcommand{\skn}{\par\smallskip\noindent}
\newcommand{\ges}{\geqslant}
\newcommand{\les}{\leqslant}
\newcommand{\one}{\hskip1pt}

\newcommand{\msum}{\hbox{$\sum$}}

\newcommand{\mprod}{\hbox{$\prod$}}
\newcommand{\CCb}{{\mathcal C}^{\ssb}}
\newcommand{\Fb}{{\mathcal F}^{\ssb}}

\newcommand{\Ec}{{\mathcal E}}
\newcommand{\Fc}{{\mathcal F}}
\newcommand{\Gc}{{\mathcal G}}
\newcommand{\M}{{\mathcal M}}
\newcommand{\OO}{{\mathcal O}}
\newcommand{\Hc}{{\mathcal H}}

\newcommand{\Sc}{{\mathcal S}}

\newcommand{\PP}{{\mathbb P}}
\newcommand{\Q}{{\mathbb Q}}
\newcommand{\C}{{\mathbb C}}
\newcommand{\DD}{{\mathbb D}}

\newcommand{\RR}{{\mathbf R}}
\newcommand{\Z}{{\mathbb Z}}

\newcommand{\pb}{{\mathbf p}}
\newcommand{\Ht}{\widetilde{H}}

\newcommand{\Ut}{\widetilde{U}}
\newcommand{\Vt}{\widetilde{V}}
\newcommand{\Xt}{\widetilde{X}}

\newcommand{\CDiv}{{\rm CDiv}}
\newcommand{\CH}{{\rm CH}}
\newcommand{\Div}{{\rm Div}}
\newcommand{\Pic}{{\rm Pic}}
\newcommand{\IC}{{\rm IC}}
\newcommand{\IH}{{\rm IH}}
\newcommand{\Gr}{{\rm Gr}}
\newcommand{\id}{{\rm id}}
\newcommand{\Sing}{{\rm Sing}}

\newcommand{\al}{\alpha}
\newcommand{\be}{\beta}
\newcommand{\ga}{\gamma}
\newcommand{\Ga}{\Gamma}

\newcommand{\la}{\lambda}

\newcommand{\si}{\sigma}
\newcommand{\Si}{\Sigma}
\newcommand{\om}{\omega}

\newcommand{\pp}{{}^{\bf p}\!}
\newcommand{\ppp}{{}^{{\bf p}^+}\!\!}
\newcommand{\pppo}{{}^{{\bf p}^+}\!\one}
\newcommand{\dd}{\partial}
\newcommand{\ddd}{{\rm d}}
\newcommand{\tos}{\,{\to}\,}
\newcommand{\eq}{\,{=}\,}
\newcommand{\defs}{\,{:=}\,}
\newcommand{\nes}{\,{\ne}\,}
\newcommand{\ins}{\,{\in}\,}
\newcommand{\sst}{\,{\subset}\,}
\newcommand{\stm}{\,{\setminus}\,}
\newcommand{\gess}{\,{\ges}\,}
\newcommand{\less}{\,{\les}\,}
\newcommand{\sgt}{\,{>}\,}
\newcommand{\slt}{\,{<}\,}
\newcommand{\col}{\,{:}\,}
\newcommand{\pl}{\one {+}\one}
\newcommand{\mi}{\one {-}\one}
\newcommand{\bl}{\bigl}
\newcommand{\br}{\bigr}

\newcommand{\ssb}{\raise.15ex\h{${\scriptscriptstyle\bullet}$}}
\newcommand{\ssc}{\,\raise.15ex\h{${\scriptstyle\circ}$}\,}

\newcommand{\into}{\hookrightarrow}
\newcommand{\simto}{\,\,\rlap{\hskip1.5mm\raise1.4mm\hbox{$\sim$}}\hbox{$\longrightarrow$}\,\,}

\makeatletter
\renewcommand\section{\@startsection{section}{1}{0pt}{-3ex plus -1ex minus -.2ex}{2.3ex plus.2ex}{\centering\normalfont\bfseries}}
\makeatother
\theoremstyle{plain}
\newtheorem{thm}{Theorem}[section]

\newtheorem{prop}[thm]{Proposition}
\newtheorem{lem}[thm]{Lemma}
\newtheorem{ithm}{Theorem}
\newtheorem{icor}{Corollary}
\newtheorem{iprop}{Proposition}

\theoremstyle{definition}
\newtheorem{rem}[thm]{Remark}

\newtheorem*{conv}{Convention}

\begin{document}
\title{Factoriality of normal projective varieties}
\author[S.-J. Jung]{Seung-Jo Jung}
\address{S.-J. Jung : Department of Mathematics Education, and Institute of Pure and Applied Mathematics, Jeonbuk National University, Jeonju, 54896, Korea}
\email{seungjo@jbnu.ac.kr}
\author[M. Saito]{Morihiko Saito}
\address{M. Saito : RIMS Kyoto University, Kyoto 606-8502 Japan}
\email{msaito@kurims.kyoto-u.ac.jp}
\thanks{This work was partially supported by National Research Foundation of Korea NRF-2021R1C1C1004097.}
\begin{abstract} For a normal projective variety $X$, the $\mathbb Q$-factoriality defect $\sigma(X)$ is defined to be the rank of the quotient of the group of Weil divisors by the subgroup of Cartier ones. We prove an improvement of a topological formula of S.G. Park and M. Popa asserting that $\sigma(X)\le h^{2n-2}(X)-h^2(X)$ by assuming only 1-semi-rational singularities instead of rational singularities, and the equality holds in the 2-semi-rational case. Here the singularities are called $j$-semi-rational if $R^k\pi_*{\mathcal O}_{\widetilde{X}}=0$ for any $k\in[1,j]$ with $\pi:\widetilde{X}\to X$ a desingularization, $h^k(X):=\dim H^k(X,{\mathbb Q})$, and $n:=\dim X\ge2$. We also show (a slight generalization of) the assertion that $\mathbb Q$-factoriality implies factoriality if $X$ is a local complete intersection whose singular locus has at least codimension three. We then get a new proof for the projective case of Grothendieck's theorem asserting that $X$ is factorial if it is a local complete intersection whose singular locus has at least codimension four. We also show that a local complete intersection $X$ of dimension 3 having only isolated singularities is factorial if the (topological) defect ${\rm def}(X):=h^4(X)-h^2(X)$ vanishes, without any assumption on rational singularities.
\end{abstract}
\maketitle

\section*{Introduction} \label{intr}
\nin
Let $X$ be a connected normal projective variety of dimension $n\gess2$. We denote by $\Div(X)$ the free abelian group of Weil divisors on $X$, which is freely generated by irreducible reduced divisors. Let $\CDiv(X)$ be the {\it subgroup\one} consisting of Cartier divisors, that is, Zariski-locally principal Weil divisors, on $X$ (since $X$ is normal). When we consider the divisors modulo linear (or rational) equivalence, that is, if we divide these abelian groups by the image of the divisor map from the rational function field $\C(X)$ of $X$, they are denoted by $\Div(X)_{\sim\rm lin}$, $\CDiv(X)_{\sim\rm lin}$, and are identified with the divisor class group ${\rm Cl}(X)$ and the Picard group $\Pic(X)\eq H^1(X,\OO^*_X)$ respectively. We will denote the rank of $\Div(X)/\CDiv(X)$ by $\si(X)$ if it is finite. It is sometimes called the {\it $\Q$-factoriality defect.}
\sk
It is known that the first Chow group $\CH^1(X)\,({=}\,{\rm Cl}(X))$ can be described in a Hodge-theoretical way, see for instance \cite[Theorem 0.3]{bconj}. We can use a more simplified argument quite effectively for the proof of the main theorem.
\sk
Set $\Si\defs\Sing\,X$. Let $\pi\col\Xt\tos X$ be a resolution of singularities such that $E\defs\pi^{-1}(\Si)$ is a divisor with simple normal crossings. In the case $X$ has only {\it rational singularities,} we have the isomorphism $\RR\pi_*\OO_{\Xt}\eq\OO_X$. This condition is however too strong. In this paper we consider the following {\rm $j$-semi-rationality\one} condition\,:
\begin{equation} \label{1}
R^k\pi_*\OO_{\Xt}=0\q\h{for any}\,\,\,k\ins\{1,\dots,j\}.
\end{equation}
We will assume this condition for $j\eq1$ or $2$ in the main theorem. We have the following.

\begin{iprop} \label{P1}
Assume $X$ is Cohen-Macaulay, for instance, a local complete intersection. Then the singularities of $X$ are $j$-semi-rational if ${\rm codim}_X{\rm NRSing}\,X\gess j{+}2$.
\end{iprop}

Here ${\rm NRSing}\,X$ denotes the {\it non-rational singularity locus of $X$}, that is, the complement of the largest Zariski-open subset of $X$ having only rational singularities (or the union of the supports of $R^k\pi_*\OO_{\Xt}$ for $k\sgt0$). This assertion can be reduced to Remark~\ref{R2.7} below using the Grauert-Riemenschneider vanishing theorem \cite{GrRi} and Grothendieck duality, since they imply that the direct image complex $\RR\pi_*\om_{\Xt}$ is quasi-isomorphic to an $\OO_X$-submodule of the dualizing sheaf $\om_X$. It is well known that in the isolated hypersurface singularity case, $X$ has an $(n{-}2)$-semi-rational singularity so that the geometric genus is defined by the dimension of $R^{n-1}\pi_*\OO_{\Xt}$, see for instance \cite{exp}.
\sk
Applying the ``functoriality" of cubic hyperresolution as in Proposition~\ref{P1.1}, we can prove an improvement of a {\it topological\one} formula of S.G.~Park and M.~Popa \cite[Theorem A]{PaPo} for $\Q$-factoriality defect (assuming rational singularities of $X$ and generalizing \cite{NaSt} in the case $n\eq3$) as follows.

\begin{ithm} \label{T1}
Let $X$ be a connected normal projective variety of dimension $n\gess 2$. If the singularities of $X$ are $1$-semi-rational as in \eqref{1}, we have the inequality
\begin{equation}  \label{2}
\si(X)\less h^{2n-2}(X)\mi h^2(X),
\end{equation}
in particular, $X$ is $\Q$-factorial in the case $h^{2n-2}(X)\eq h^2(X)$ $($for instance if $n\eq2)$.
\sk
If the singularities of $X$ are $2$-semi-rational, we have the equality
\begin{equation} \label{3}
\si(X)\eq h^{2n-2}(X)\mi h^2(X).
\end{equation}
\end{ithm}

Here $h^k(X)\defs\dim_{\Q}H^k(X,\Q)$ for $k\ins\Z$. Recall that $X$ is factorial (resp. $\Q$-factorial) if and only if $\Div(X)\eq\CDiv(X)$ (resp. $\si(X)\eq0$). We can show that these are equivalent if $X$ is a local complete intersection and ${\rm codim}_X\Sing\,X\gess3$, see Theorem~\ref{T4.1} below. So $X$ is factorial if $X$ is locally a hypersurface having only isolated singularities with $n\eq3$ and 1 is not an eigenvalue of the Milnor monodromy (so that $\Q_X[n]\simto{\rm IC}_X\Q$ and $\si(X)\eq0$, see also \cite[Theorem 1]{JS}) even if the singularities of $X$ are {\it not\one} 2-semi-rational, see Remark~\ref{R2.10} below. Using Theorems~\ref{T1} and \ref{T4.1} below with Proposition~\ref{P1}, we can get the following (which is a slight improvement of Grothendieck's theorem in the projective case, see \cite[XI, Theorem 3.13 (ii)]{Gro}, where one assumes ${\rm codim}_X \Sing\,X\gess4$, see also Remark~\ref{R4.3} below).

\begin{icor} \label{C1}
A normal projective variety $X$ is factorial if it is a local complete intersection with ${\rm codim}_X \Sing\,X\gess3$ and if $h^{2n-2}(X)\eq h^2(X)$ $($for instance if ${\rm codim}_X{\rm HSing}\,X\gess4)$.
\end{icor}

Here ${\rm HSing}\,X$ is the {\it $\Q$-homology singular locus,} that is, the complement of the largest open subvariety of $X$ which is a $\Q$-homology manifold. Note that a variety $X$ is a $\Q$-homology manifold if and only if we have the canonical isomorphism $\Q_X[n]\simto\IC_X\Q$. (This follows easily from \cite[(4.5.6--9)]{mhm} using the self-duality of $\IC_X\Q_h$.)
The right-hand side of \eqref{2} then vanishes under the last assumption of Corollary~\ref{C1} using Lemma~\ref{L2.1} below, see also \cite{JS} for $n\eq3$. There are many examples where the hypotheses of Corollary~\ref{C1} are satisfied with $n\eq3$, see Remark~\ref{R2.10} below. By a similar argument we can get the following.

\begin{icor} \label{C2}
A normal projective variety $X$ is $\Q$-factorial in the case $X$ is Cohen-Macaulay, ${\rm codim}_X{\rm NRSing}\,X\gess3$, ${\rm codim}_X{\rm HSing}\,X\gess4$, and $\pp\Hc^k(\Q_X[\dim X])\eq0$ for any $k\nes0$.
\end{icor}

Here $\pp\Hc^{\ssb}$ is the cohomology functor in \cite{BBD}. Note that $\Q$-factorial cannot be replaced with factorial, considering the case of the cone of a Veronese embedding, see Remark~\ref{R2.8} below (with $n\eq4$).
Using Theorems~\ref{T1} and \ref{T4.1} below, we can deduce also the following.

\begin{ithm} \label{T2}
Assume $X\sst\PP^N$ is a projective hypersurface, or more generally, a $($global$\one)$ complete intersection, such that ${\rm codim}_X\Sing\,X\gess3$, ${\rm codim}_X{\rm NRSing}\,X\gess4$, $n\eq\dim X\gess3$, and $X$ is not contained in any hyperplane of $\PP^N$ $($replacing $\PP^N$ if necessary$)$. For a hyperplane $H\sst\PP^N$, let $H_X$ be the restriction of $H$ to $X$ as a Cartier divisor. The following four conditions are equivalent$\,:$
\skn
{\rm(a)} The variety $X$ is factorial.
\skn
{\rm(b)} The equality $h^{2n-2}(X)\eq h^2(X)$ holds.
\skn
{\rm(c)} The divisor class group ${\rm Cl}(X)\eq\Div(X)_{\sim{\rm lin}}$ is freely generated by the class of $H_X$ for some $($or any$)$ hyperplane $H\sst\PP^N$.
\skn
{\rm(d)} The affine ring $A_H$ of the affine variety $X\stm H_X$ is a unique factorization domain for some $($or any$)$ hyperplane $H\sst\PP^N$.
\sk
In the case $X$ is a projective cone of a $($global$\one)$ complete intersection $Y\sst\PP^{N-1}$ with $\dim Y\gess3$ $($where ${\rm codim}_Y\Sing\,Y\gess3$, ${\rm codim}_Y{\rm NRSing}\,Y\gess4)$, we have the equality
\begin{equation} \label{4}
h^{2n-2}(X)\mi h^2(X)\eq h^{2n-4}(Y)\mi h^2(Y),
\end{equation}
so that any of the above equivalent four conditions for $X$ is further equivalent to any of those conditions for $Y$ with $n,N$ replaced by $n{-}1,N{-}1$ respectively.
\end{ithm}

This generalizes some assertion in \cite[Remark 1.2]{Che} for hypersurfaces $Y$ with $\dim Y\eq3$, see also \cite[Exer.\ II.6.3]{Ha}, \cite{Ho}.
In conditions (c), (d), ``some" may be replaced by ``any", that is, the assertion holds for any $H$ if it holds for some $H$. (This is clear for (c).) The equivalence between (c) and (d) holds for any reduced projective variety with $H$ fixed. In order that condition (c) is satisfied in the $X$ smooth case, the Picard number of $X$ must be one and $h^1(X)$ must vanish, although these are not sufficient (considering a Veronese case, where a hyperplane section is a hypersurface of degree $m$, see also Remark~\ref{R2.8} below) if we do not assume the condition that $X$ is a (global) complete intersection.
\sk
The proof of Theorem~\ref{T2} uses Artin's vanishing theorem and semi-perversity with integral coefficients, see Remark~\ref{R3.1} below. It is interesting that rectified homological depth, which is very closely related to semi-perversity (see \cite{HaLe}), was explored in \cite[XIII.6]{Gro} after the study of purity and factoriality in \cite[X.3.4 and XI.3.13]{Gro}.
\sk
We thank M.~Popa, J.~Koll\'ar, and A.~Petracci for informing us of Grothendieck's work on factoriality \cite[XI.3.13]{Gro}, the one on purity \cite[X.3.4]{Gro}, and an absence of some additional assumption in a previous version of Corollary~\ref{C2} respectively. 

\begin{conv}
An algebraic variety means a separated reduced scheme of finite type over $\C$ where only the closed points are considered, that is, a variety in the sense of J.-P.~Serre. However we use also the classical topology to employ constructible sheaves and analytic sheaves especially on non-algebraic open subsets. This does not cause a problem by applying GAGA in the projective case.
\end{conv}

\tableofcontents
\numberwithin{equation}{section}

\section{Restriction on Hodge types} \label{S1}

We denote by $\Q_{h,X}$ the pullback $a_X^*\Q_h\ins D^b{\rm MHM}(X)$, where $a_X\col X\tos{\rm pt}$ is the structure morphism and $\Q_h$ is the mixed Hodge structure of rank~1 and weight~0, which is identified with a mixed Hodge module on $\rm pt$, see \cite{mhm}. Let $\IC_X\Q_h$ be the intersection complex Hodge module. In this and next section we set  
\begin{equation*}
\CCb:=C\bl(\Q_{h,X}[n]\to\IC_X\Q_h\br),
\end{equation*}
(which is different from the one in Section~\ref{S4} below). Put
\begin{equation*}
L^k\defs\Hc^k\bl(\IC_X\Q[-n]\br)\eq\Hc^k\bl({\rm rat}(\CCb)[-n]\br)\q\h{for}\,\,\, k\ins[1,n{-}1],
\end{equation*}
and $L^k\eq0$ otherwise, where ${\rm rat}(\CCb)$ denotes the underlying $\Q$-complex. The stalks $L^k_x$ can be obtained by using $H^ki_x^*$ with $i_x\col\{x\}\into X$ the inclusion, and they have weights at most $k$ for $x\ins\Si$, see \cite[(4.5.2) and also Remark 4.6.1]{mhm}.
We have the following.

\begin{prop} \label{P1.1}
Assume the singularities of $X$ are $j$-semi-rational $($see \eqref{1}$)$, where $j\eq1$ or $2$. Then for any $k\less j$, we have $\Gr^0_FL^k_x\eq0$ $($in particular, $L^1_x\eq0)$ and isomorphisms
\begin{equation*}
H^k(X,\OO_X)\simto\Gr^0_FH^k(X)_{\C}\simto\Gr^0_FH^k(\Xt)_{\C}=\Gr^0_F\IH^k(X)_{\C},
\end{equation*}
where the first isomorphism holds also for $k\eq j{+}1$.
\end{prop}

\begin{proof}
We set $H^k(X)_A\defs H^k(X,A)$ for a ring $A$ in general. For $k\ins[1,j{+}1]$, we have the commutative diagram
\begin{equation} \label{1.1}
\begin{gathered}
\xymatrix{H^k(X)_{\C} \ar[r] \ar@{->>}[rd] & H^k(X,\OO_X)\, \ar@{^{(}->}[r]^{{}^k\pi_{\OO}^*} \ar@{^{(}->>}[d] & H^k(\Xt,\OO_{\Xt}) \ar@{=}[d]\\ & \Gr^0_FH^k(X)_{\C}\, \ar@{^{(}->}[r]^{{}^k\pi_{\C}^*} & \Gr^0_FH^k(\Xt)_{\C}}
\end{gathered}
\end{equation}
using the ``functoriality" of {\it cubic hyperresolutions\one} for the morphism $\Xt\tos X$, see \cite{GNPP}. (As for the commutativity of the left part of the diagram, note that there is a morphism from $\OO_X$ to the derived direct image of the structure sheaves of the cubic resolution in a compatible way with a morphism from $\C_X$ to the derived direct image of the constant sheaves.)
We see that $^k\pi_{\OO}^*$ is {\it bijective\one} for $k\ins[1,j]$, and {\it injective\one} for $k\eq j{+}1$, assuming the $j$-semi-rationality (see \eqref{1}) and using the Leray-type spectral sequence
\begin{equation} \label{1.2}
E^{p,q}_2=H^p(X,R^q\pi_*\OO_{\Xt})\Longrightarrow H^{p+q}(\Xt,\OO_{\Xt}).
\end{equation}
(Here we do not employ the assertion that rational singularities are du Bois, see \cite{Ko}, \cite[Theorem 5.4]{mhc}.)
\sk
The last isomorphisms of Proposition~\ref{P1.1} then follow using the canonical isomorphisms
\begin{equation} \label{1.3}
\Gr^0_FH^k(\Xt)_{\C}=\Gr^0_F\IH^k(X)_{\C}\q\h{for any}\,\,\,k\ins\Z.
\end{equation}
Indeed, as a consequence of the decomposition theorem for Hodge modules (see \cite[(4.5.4)]{mhm}), there are non-canonical isomorphisms of mixed Hodge structures
\begin{equation*}
H^k(\Xt)_{\Q}\cong H_1^k\oplus H_2^k\q\h{for any}\,\,\,k\ins\Z,
\end{equation*}
with $H_1^k\defs\IH^k(X)_{\Q}$. We have moreover
\begin{equation} \label{1.4}
\Gr^p_FH_2^k\eq0\q\h{for}\,\,\,p=0,n\,\,\,\h{and any}\,\,\,k\ins\Z.
\end{equation}
Indeed, the assertion for $p\eq0$ is reduced to the case $p\eq n$ by the self--duality of $H^{\ssb}_2$, and the latter case follows from Koll\'ar's torsion-freeness in the Hodge module theory (see for instance \cite{kol}), since the direct factors of $\pi_*\Q_{h,\Xt}[n]$ are supported on $\Si\sst X$ except $\IC_X\Q_h$. So the canonical isomorphism \eqref{1.3} follows using the filtration induced by the truncations $\tau_{\les\ssb}$ on $\pi_*\Q_{h,\Xt}[n]$ (since the {\it strict support decomposition\one} of pure Hodge modules is canonical).
\sk
We now show that $\Gr^0_FL^k_x\eq0$ for any $k\less j$ by induction on $j$. Let $Z$ be the set of $x\ins\Si$ with $\Gr^0_FL^j_x\nes0$. It contains a Zariski-locally closed subset $Z'$ of $\Si$ whose closure contains $Z$. (Note that the restrictions of the cohomology sheaves $\Hc^k\IC_X\Q$ to each stratum of a Whitney stratification underlie variations of mixed Hodge structure using \cite[Remark 4.6.2]{mhm}.) We can verify this by applying the decomposition theorem to $\pi'_*\Gr^W_k\Q_{h,E}$ with $\pi'\defs\pi|_E$. Taking a sufficiently general hyperplane section, we may assume $\dim Z'\eq0$ if $Z\nes\emptyset$. The assertion for $L^k_x$ for $k\eq j$ then follows from the exact sequence of mixed Hodge structures
\begin{equation} \label{1.5}
H^k(X)_{\Q}\to\IH^k(X)_{\Q}\buildrel{\ga^k}\over\to\Ga(\Si,L^k)\to H^{k+1}(X)_{\Q}\to\IH^{k+1}(X)_{\Q},
\end{equation}
using the {\it inductive hypothesis\one} in the case $j\eq2$. Indeed, these imply that $\Gr_F^0L^k_x\eq0$ in view of \eqref{1.3} and the bijectivity (resp. injectivity) of $^k\pi^*_{\C}$ for $k\ins[1,j]$ (resp. $k\eq j{+}1$) in \eqref{1.1}. This finishes the proof of Proposition~\ref{P1.1}.
\end{proof}

\begin{rem} \label{R1.2}
Using the functorial morphism $\id\tos (i'_x)_*i_x^{\prime\,*}$ with $i'_x\col\{x\}\,{\into}\,\Si$ the inclusion for $x\ins\Si$, we can get an assertion similar to Proposition~\ref{P1.1} for $\Ga(\Si,L^2)$.
\end{rem}

\begin{rem} \label{R1.3}
Assuming 1-semi-rationality, the equality \eqref{3} is equivalent to the injectivity of the $E_3$-differential $\ddd_3^{0,2}\col E_3^{0,2}\eq E_2^{0,2}\tos E_3^{3,0}\eq E_2^{3,0}$ of the spectral sequence \eqref{1.1} employing $E_2^{p,1}\eq0$ for any $p\ins\Z$. Indeed, the kernel of $\ddd_3^{0,2}$ is canonically isomorphic to $\Gr_F^0H'_{\C}$ in the last part of the proof of Theorem~\ref{T1} in Section~\ref{S2} below, using the commutative diagram \eqref{1.1}, the spectral sequence \eqref{1.2}, and the canonical isomorphism \eqref{1.3}. Note that the image of $\ddd_3^{0,2}$ in $E_2^{3,0}$ can be nonzero and the converse of the last assertion of Theorem~\ref{T1} can be false as shown in Remark~\ref{R2.10} below.\end{rem}

\section{Proof of Theorem~\ref{T1} and some remarks} \label{S2}

We first consider the 1-semi-rationality case.
By the vanishing of $L^1$ in Proposition~\ref{P1.1}, there are isomorphisms of $\Q$-Hodge structures
\begin{equation} \label{2.1}
H^1(X)_{\Q}\eq\IH^1(X)_{\Q}\,{\cong}\,\IH^{2n-1}(X)_{\Q}(n{-}1)\eq H^{2n-1}(X)_{\Q}(n{-}1),
\end{equation}
where the middle isomorphism depends on the choice of an ample line bundle on $X$, and the last two isomorphisms hold also for $k\eq1$ assuming only the 1-semi-rationality, since $\dim{\rm supp}\,\CCb\less n{-}3$ (using only that $L^1\eq0$), see Lemma~\ref{L2.1} below. Combining this with \eqref{1.5} for $k\eq1$, we then get
\begin{equation} \label{2.2}
H^2(X)_{\Q}\into\IH^2(X)_{\Q}\,{\cong}\,\IH^{2n-2}(X)_{\Q}(n{-}2)\eq H^{2n-2}(X)_{\Q}(n{-}2).
\end{equation}
In particular, $H^k(X)_{\Q}$, $H^{2n-k}(X)_{\Q}$ are {\it polarizable Hodge structures\one} of pure weight $k$ and $2n{-}k$ respectively if $k\eq1$ or 2.
\sk
Using the long exact sequence associated with the exponential sequence together with Proposition~\ref{P1.1} and \eqref{2.2}, we can then deduce the short exact sequence
\begin{equation} \label{2.3}
\begin{aligned}
&0\to\Pic^0(X)\to\Pic(X)\to{\rm Hom}_{\rm HS}\bl(\Z,H^2(X)_{\Z}(1)\br)\to0,\\
&\h{with}\q\q\Pic^0(X)=J\bl(H^1(X)_{\Z}\br):=\Gr^0_FH^1(X)_{\C}/H^1(X)_{\Z},
\end{aligned}
\end{equation}
and $\Pic^0(X)$ is an abelian variety. Here ${\rm Hom}_{\rm HS}$ denotes the abelian group of morphisms of Hodge structures.
\sk
We have an inclusion of abelian varieties $\Pic^0(X)\into\Pic^0(\Xt)\eq\CH^1(\Xt)_{\rm alg}$, where the last group denotes the subgroup of $\CH^1(\Xt)$ consisting of divisors algebraically equivalent to zero, see \cite{Fu}. (Indeed, its composition with the pushforward $\pi_*\col\CH^1(\Xt)_{\rm alg}\tos\CH^1(X)_{\rm alg}$ coincides with the natural morphism, which is injective by construction.) This inclusion is an {\it isogeny of abelian varieties\one} (hence it is surjective considering the morphism of the universal coverings), since the source and the target have the same dimension by Proposition~\ref{P1.1}. We thus get the isomorphism
\begin{equation} \label{2.4}
\Pic^0(X)\simto\Pic^0(\Xt)=\CH^1(\Xt)_{\rm alg}.
\end{equation}
\sk
Let $E(\Xt)\sst\Div(\Xt)$ be the kernel of the pushforward $\pi_*\col\Div(\Xt)\tos\Div(X)$. This is freely generated by exceptional divisors. Its intersection with the image of $\C(\Xt)$ is zero (taking the direct image on $X$ of $g\ins\C(\Xt)$ with ${\rm div}(g)\ins E(\Xt)$ and using the normality of $X$), hence $E(\Xt)$ is identified with the kernel of the {\it surjection\one} $\pi_*\col\CH^1(\Xt)\tos\CH^1(X)$. We have moreover
\begin{equation} \label{2.5}
E(\Xt)\cap\CH^1(\Xt)_{\rm alg}=0,
\end{equation}
hence $E(\Xt)$ is identified with a subgroup of ${\rm NS}(\Xt)\eq\CH^1(\Xt)/\CH^1(\Xt)_{\rm alg}$.
Indeed, if the intersection is nonzero, it implies the non-injectivity of $\Pic^0(X)\to\CH^1(X)$ considering the composition of the isomorphism \eqref{2.4} with the restriction of $\pi_*\col\CH^1(\Xt)\tos\CH^1(X)$ to $\CH^1(\Xt)_{\rm alg}$. (This gives a simple proof of an assertion considered in Remark~\ref{R2.2} below in the 1-semi-rational singularity case.)
\sk
We then see that the injective image of $\CH^1(\Xt)_{\rm alg}$ in $\CH^1(X)$ coincides with $\CH^1(X)_{\rm alg}$ using the {\it divisibility\one} of $\CH^1(X)_{\rm alg}$ (see \cite[Example 19.1.2]{Fu}), since a quotient group of a divisible group is divisible and a nonzero finitely generated abelian group is not divisible. Actually this coincidence is not needed for the proof. Indeed, there is an exact sequence of {\it pure\one} Hodge structures
\begin{equation} \label{2.6}
0\to H^{2n-2}(X)_{\Q}\to H^{2n-2}(\Xt)_{\Q}\to H^{2n-2}(E)_{\Q}\to0,
\end{equation}
which can be proved by using the long exact sequences associated to the cohomologies with compact supports of $\Xt\stm E\eq X\stm\Si$ as well as \eqref{2.1}--\eqref{2.2}, since $H_c^{2n-2}(X\stm\Si)_{\Q}\eq H^{2n-2}(X)_{\Q}$ (because $\dim\Si\less n{-}2$), see also \cite[Proposition 3.9]{dBo}, \cite{GNPP}. Note that $H^{2n-3}(E)_{\Q}$ has weights at most $2n{-}3$, since $E$ is compact. Setting $E_{\rm sm}\defs E\stm\Sing\,E$, we have moreover the isomorphisms
\begin{equation} \label{2.7}
H^{2n-2}(E)_{\Q}\eq H^{2n-2}_c(E_{\rm sm})_{\Q}\eq E(\Xt)_{\Q}(1{-}n).
\end{equation}
\sk
By the compatibility of the cycle class map with the direct image by proper morphisms (which can be verified modulo torsion by employing the integration of $C^{\infty}$ forms on cycles), we then get the short exact sequence
\begin{equation} \label{2.8}
0\to\Pic^0(X)_{\Q}\to\CH^1(X)_{\Q}\to\bl({\rm NS}(\Xt)/E(\Xt)\br){}_{\Q}\to0,
\end{equation}
and the last term is identified with the image of $\CH^1(X)_{\Q}$ in $H_{2n-2}(X)_{\Q}(1{-}n)$ and with
\begin{equation} \label{2.9}
{\rm Hom}_{\rm HS}\bl(\Q,H_{2n-2}(X)_{\Q}(1{-}n)\br)\,\bl({\cong}\,{\rm Hom}_{\rm HS}\bl(\Q,H^{2n-2}(X)_{\Q}(n{-}1)\br)\br),
\end{equation}
using \eqref{2.6}--\eqref{2.7} and semisimplicity of polarizable $\Q$-Hodge structures, since
\begin{equation*}
{\rm NS}(\Xt)_{\Q}={\rm Hom}_{\rm HS}\bl(\Q,H_{2n-2}(\Xt)_{\Q}(1{-}n)\br),
\end{equation*}
with $M_{\Q}\defs M{\otimes}_{\Z}\Q$ for an abelian group $M$ in general. For this we show that the quotient group $\CH^1(X)_{{\rm hom},\Q}/\Pic^0(X)_{\Q}$ is not isomorphic to a {\it nonzero\one} submodule of $\bl({\rm NS}(\Xt)/E(\Xt)\br){}_{\Q}$ via \eqref{2.8} by using the {\it surjectivity\one} of the dual of the first morphism of \eqref{2.6} and counting the dimensions. Here $\CH^1(X)_{\rm hom}$ is the subgroup consisting of divisors homologically equivalent to zero, which coincides with the kernel of the cycle class map to $H_{2n-2}(X)_{\Z}(1{-}n)$. (This argument implies the coincidence of $\CH^1(X)_{\rm alg}$ and $\CH^1(X)_{\rm hom}$ modulo torsion, that is, after taking the tensor product with $\Q$ over $\Z$, see for instance \cite[remark after Theorem 2.1]{PaPo} and also \cite{bconj}.) 
\sk
The assertion for $j\eq1$ now follows comparing \eqref{2.3} with \eqref{2.8}--\eqref{2.9} and using \eqref{2.2}. Here we get the inequality, since $\dim_{\Q}{\rm Hom}_{\rm HS}(\Q,H'_{\Q}(1))\less\dim_{\Q}H'_{\Q}$, where $H'_{\Q}$ denotes the cokernel of the first inclusion in \eqref{2.2}.
In the case $j\eq2$, the cokernel $H'_{\Q}$ has type $(1,1)$ in view of Proposition~\ref{P1.1}, hence we get the equality using semisimplicity of polarizable $\Q$-Hodge structures. This finishes the proof of Theorem~\ref{T1}.
\sk
We note here the following well-known assertion (which is shown by using the vanishing of $H^k(X,\M')$ for any $k\notin[-m,m]$, where $\M'$ is a mixed Hodge module with $\dim{\rm supp}\,\M'\eq m$).

\begin{lem} \label{L2.1}
Let $r$ be a non-negative integer with $\pp\Hc^k(\Q_X[n])\eq0$ for any $k\,{\notin}\,[-r,0]$ $($for instance, $X$ is a local complete intersection with $r\eq0)$. Set $c\defs{\rm codim}_X{\rm HSing}\,X$. Then $H^k\CCb\eq0$ for any $k\,{\notin}\,[-r{-}1,-1]$, and $H^k(\Si,\CCb)\eq0$ for any $k\,{\notin}\,[-n{+}c-r{-}1,n{-}c{-}1]$.
\end{lem}

\begin{rem} \label{R2.2}
Proof of Theorem~\ref{T1} in \cite{PaPo} uses linear independence of the topological cycle classes of exceptional divisors of a desingularization of $X$ (among others). This can be reduced to negative definiteness of the intersection matrix of exceptional divisors (see \cite[p.~6]{Mu} where only an outline of proof is indicated; a converse assertion is shown in \cite[p.~367]{Gra}) by taking a general hyperplane section of $\Xt$ or the pullback of one for $X$ depending on the coefficients of a linear dependence. In the first case it seems rather difficult to proceed by induction (finding an appropriate blow-down); it might be simpler to consider an iterated hyperplane section of dimension~2 together with the normalization of its image in $X$ forgetting the iterated hyperplane sections of exceptional divisors whose images in $X$ are not 0-dimensional.
\end{rem}

\begin{rem} \label{R2.3}
There are many examples with $L^1\nes0$ and $h^1(X)\nes h^{2n-1}(X)$ (see \cite[Theorem A]{PaPo}) even in the isolated singularity case. For instance, assume $X$ is the projective cone of a smooth projective variety $Y$ with vertex $0$ and $H^1(Y)\nes0$. Let $\Xt$ be the blow-up of $X$ along 0, which is a $\PP^1$-bundle over $Y$. We then get
\begin{equation} \label{2.10}
\begin{aligned}
H^k(\Xt)&=H^k(Y)\oplus H^{k-2}(Y)(-1),\\
H^k(X)&=\begin{cases} H^{k-2}(Y)(-1)&\h{if}\,\,\,k\nes0,\\ \Q&\h{if}\,\,\,k\eq0,
\end{cases}
\end{aligned}
\end{equation}
using the long exact sequence $H^k_c(X\stm\{0\})\to H^k(X)\to H^k(\{0\})\to$. (These isomorphisms hold without assuming the smoothness of $Y$.) Here the $L^k_0$ are given by the primitive part of $H^k(Y)$, and the non-primitive part gives the direct factor of the direct image $\pi_*\Q_{h,\Xt}[n]$ supported on 0. This follows from the Thom-Gysin sequence, see for instance \cite{RSW}. We then see that $h^1(X)\nes h^{2n-1}(X)$ if $h^1(Y)\nes0$, see \cite[Theorem A]{PaPo}.
\end{rem}

\begin{rem} \label{R2.4}
With the notation of Remark~\ref{R2.3}, assume $Y$ is an elliptic curve $E\sst\PP^2$ defined by $x^3\pl y^3\pl z^3\eq0$ for instance. We see that the divisor class $[D_e]\mi[D_{e_0}]$ is not a torsion and has infinite order if $e\ins E$ is not a torsion point (restricting to the zero section of the line bundle $X\stm\{0\}\tos E$). Here $e_0\ins E$ is the origin (which intersects a line in $\PP^2$ with multiplicity three; for instance $e_0\eq[0:-1:1]$ with line defined by $y\pl z\eq0$) and $D_e,D_{e_0}$ are the closures of the inverse images of $e,e_0$ by the projection $X\stm\{0\}\tos E$. Recall that the sum of the three intersection points of $E$ and a line in $\PP^2$ is zero by the rule of addition of elliptic curves. Here $X$ is normal using Serre's condition, see \cite[p.~185]{Ha}. (In the non-hypersurface case it does not seem clear whether the cone $X$ is normal; it may be necessary to take the normalization which does not change the underlying topological space.) This example can be generalized to the case of curves of higher genus or smooth projective varieties $Y$ with $H^1(Y,\OO_Y)\nes0$.
\end{rem}

\begin{rem} \label{R2.5}
In the notation of Remark~\ref{R2.3}, assume $Y\sst\PP^3$ is a smooth surface such that $d\defs\deg Y\gess4$, that is, $H^2(Y,\OO_Y)\nes0$. We have $\si(X)\eq\rho(Y)\mi1$ where $\rho(Y)$ is the Picard number, that is, $\dim_{\Q}H^2(Y,\Q)\cap F^1H^2(Y,\C)$, hence $\si(X)\less h^4(X)\mi h^2(X)\mi2\one h^{2,0}(Y)$ (in particular $h^4(X)\sgt h^2(X)$) using the Hodge decomposition (and symmetry) of $H^2(Y,\C)$, where $h^{2,0}(Y)\eq\dim_{\C}H^2(Y,\OO_Y)$. (Note that $\rho(Y)\eq1$ if $Y$ is very general by the Noether-Lefschetz theorem, see for instance \cite{BN}.) Indeed, the Picard group $\Pic(X\stm\{0\})$ is isomorphic to $\Pic(Y)$ via the pullbacks by the projection and the zero-section of the line bundle, see for instance \cite[Proposition 1.9]{Fu}. Moreover an element of $\Pic(X\stm\{0\})$ is extendable to $X$ if and only if it is a multiple of a hyperplane section class (since the rational function fields of $X$ and $X\stm\{0\}$ are the same) using the surjectivity of $\OO_{\C^4,0}\tos\OO_{X,0}$ (see Remark~\ref{R3.1} below) and the Hartogs-type lemma.
\sk
By the assumption $d\gess4$, the singularity of the cone $X$ is {\it not\one} 2-semi-rational, and these give examples where the equality \eqref{2} in Theorem~\ref{T1} does not hold without this 2-semi-rationality hypothesis for $n\eq3$.
\end{rem}

\begin{rem} \label{R2.6}
There is an example of a normal projective hypersurface $X$ having non-1-semi-rational singularities with $h^1(X)\eq h^{2n-1}(X)$ (see \cite[Theorem A]{PaPo}) and $n\eq3$. We can consider for instance the case $X\sst\PP^4$ is the projective cone of a surface $Y\sst\PP^3$ defined by $x^2yz^2\pl xy^4\pl w^5\eq0$. Here $X$ has {\it non-isolated\one} singularities, and $X,Y$ are $\Q$-homology manifolds. The normality of $X$ follows from Serre's condition as explained in Remark~\ref{R2.4}. The curve $Z\sst\PP^2$ defined by $x^2yz^2\pl xy^4\eq0$ has two singular points whose spectral numbers are $j/7$ with $j\ins\{4,\dots,10\}$ and $j/8$ with $j\ins\{5,7,8,9,11\}$, see for instance \cite[Section A1]{wh} or \cite{Sing}. These do not produce an integer by adding $i/5$ for any $i\ins\{1,\dots,4\}$ to any of them. So we can apply the Thom-Sebastiani type theorem for spectrum (see for instance \cite{SS}) to show that $Y$ is a $\Q$-homology manifold. The latter theorem implies also the non-1-semi-rationality of $\Sing\,Y$ using \cite{exp} (since $\dim Y\eq2$). For the calculation of the stalk of the intersection complex at the vertex of $X$, we can employ the argument at the end of Remark~\ref{R2.3}.
\end{rem}

\begin{rem} \label{R2.7}
It seems well known that for a coherent sheaf $\Fc$ on a Cohen-Macaulay variety $Z$, we have $\Ec xt^j_{\OO_Z}(\Fc,\om_Z)\eq0$ for any $j\slt c\defs{\rm codim}_Z\one{\rm supp}\,\Fc$. This follows for instance from \cite[Lemma III.7.3]{Ha}. Recall that the dual $\DD\one\Fc$ is defined by $\RR\Hc om_{\OO_Z}(\Fc,\om_Z[\dim Z])$. This is independent of the choice of a Cohen-Macaulay variety $Z$ by Grothendieck duality for closed immersions. Hence we may assume $Z$ is {\it smooth.}
We note here a simple argument which works also in the analytic case for the convenience of the reader.
\sk
Since the assertion is local, we may assume $\Fc\eq\OO_Y$ with $Y$ a closed subvariety (using a filtration associated with a local finite presentation of $\Fc$). We may assume further that $Y$ is reduced using some filtration (to ensure the surjectivity of the last morphism of the exact sequence below). There is locally a complete intersection $Y'$ (not necessarily reduced) of the same dimension as $Y$ and containing $Y$ so that we have locally a short exact sequence
\begin{equation*}
0\to\Gc\to\OO_{Y'}\to\Fc\to0.
\end{equation*}
Since $\Ec xt^j_{\OO_Z}(\OO_{Y'},\om_Z)\eq0$ for any $j\nes c$, we get the isomorphisms
\begin{equation*}
\Ec xt^{j-1}_{\OO_Z}(\Gc,\om_Z)=\Ec xt^j_{\OO_Z}(\Fc,\om_Z)\q\h{for any}\,\,\,j\slt c.
\end{equation*}
We can then prove by decreasing induction on $c\less\dim Z$ and increasing induction on $b\slt c$ that $\Ec xt^j_{\OO_Z}(\Fc,\om_Z)\eq0$ for any $j\less b$ and any coherent sheaf $\Fc$ with ${\rm codim}_Z\one{\rm supp}\,\Fc\gess c$ (where $\Fc$ is not fixed).
\end{rem}

\begin{rem} \label{R2.8}
Let $v_m:\PP^{n-1}\into\PP^{N-1}$ be the Veronese embedding of degree $m\gess2$, where $N\eq\binom{n+m-1}{m}$ and $n\gess2$. Let $Y\sst\PP^{N-1}$ be the image of $v_m$. Let $X\sst\PP^N$ be the projective cone of $Y$. Then $X$ is Cohen-Macaulay, and Gorenstein if $n/m\ins\Z$. Indeed the affine cone of $Y$ is the quotient of $\C^n$ by an action of $\mu_m\defs\{\la\ins\C^*\mid\la^m\eq1\}\cong\Z/m\one\Z$ via the scalar multiplication on $\C^n$ (since any monomial of degree divisible by $m$ is a product of monomials of degree $m$). Hence $X$ has only rational singularities. Using the Grauert-Riemenschneider vanishing theorem and Grothendieck duality, this implies that $X$ is Cohen-Macaulay, and Gorenstein in the case $n/m\ins\Z$ (considering the action on $n$-forms). Note that the blowup of the affine cone of $Y$ at the vertex is a line bundle over $Y$. This is the quotient by $\mu_m$ of the blowup of $\C^n$ at 0 (which is a line bundle corresponding to the invertible sheaf $\OO_{\PP^{n-1}}(-1)$) and corresponds to $\OO_{\PP^{n-1}}(-m)$. Here the quotient group $\Div(X)/\CDiv(X)$ is isomorphic to $\Z/m\one\Z$. Indeed, passing to the linear equivalent classes, the divisor class group ${\rm Cl}(X)\eq\Div(X)_{\sim{\rm lin}}$ is generated by the class of the cone $\Ht$ of a hyperplane $H$ in $\PP^{n-1}\eq Y$, but the Picard group $\Pic(X)\eq\CDiv(X)_{\sim{\rm lin}}$ by the class of a hyperplane section of $X$ in $\PP^N$ or that of the cone of a hyperplane section of $Y$ in $\PP^{N-1}$, which is equal to $m[\Ht]$ (using an argument similar to the one in the end of Remark~\ref{R2.5}). Note also that the fundamental group of the link of the vertex of the projective cone $X$ is isomorphic to $\Z/m\one\Z$, and $X$ is a $\Q$-homology manifold (considering the decomposition of the direct image of $\Q_{\C^n}[n]$ by the covering transformation group $\mu_m$). The cone $X$ is {\it not\one} a local complete intersection in the case $n\gess3$ by \cite[Theorem A]{KaWa} (where \cite[X.3.3--4]{Gro} is used in the proof). This assertion follows also from Theorem~\ref{T4.1} below.
\end{rem}

\begin{rem} \label{R2.9}
For $x\ins\Si$, let $U_x$ be a sufficiently small {\it contractible\one} Stein open neighborhood (by intersecting $X$ with a sufficiently small ball $B_x$ in an ambient smooth space). Set $\Ut_x\defs\pi^{-1}(U_x)$. We have the following isomorphism which generalizes \cite[(3.1)]{NaSt} and expresses the quotient group $\Div(X)/\CDiv(X)$ via the Picard group of a desingularization and the local analytic Picard groups at singular points.
\begin{equation*}
\begin{aligned}
&\Div(X)/\CDiv(X)\\&\simto{\rm Im}\bl(\Pic(\Xt)\buildrel{\!\al}\over\to\mprod_{x\in\Si}\,\Pic(\Ut_x)/\bl(\msum_i\,\Z[E_i\cap\Ut_x]\br)\br).
\end{aligned}
\end{equation*}
Here $\Pic(\Ut_x)\defs H^1(\Ut_x,\OO^*_{\Ut_x})$, and $[E_i\cap\Ut_x]$ denotes the class of the line bundle associated with the divisor $E_i\cap\Ut_x\sst\Ut_x$(corresponding to the invertible sheaf $\OO_{\Ut_x}(E_i\cap\Ut_x))$.
\sk
If the $2$-semi-rationality condition \eqref{1} is satisfied, we get the equality
\begin{equation*}
\begin{aligned}
&\Div(X)/\CDiv(X)\\
&\simto{\rm Im}\bl(H^{1,1}(\Xt,\Z)\buildrel{\!\be}\over\to\Ga\bl(\Si,R^2\pi'_*\Z_E/\bl(\msum_i\,\Z[E_i]_{\pi'}\br)\br)\br).
\end{aligned}
\end{equation*}
Here $H^{1,1}(\Xt,\Z)\defs{\rm Ker}\bl(H^2(\Xt,\Z)\tos H^2(\Xt,\OO_{\Xt})\br)$, $\pi'\defs\pi|_E$, and $[E_i]_{\pi'}$ denotes the section of $R^2\pi'_*\Z_E$ defined by the restrictions of the Chern class of $[E_i]$ to $\pi^{-1}(x)$ for $x\ins\Si$.
Indeed, using the long exact sequence associated with the exponential sequence, we can get the isomorphisms
\begin{equation*}
\Pic(\Ut_x)=H^2(\Ut_x,\Z)=H^2(\pi^{-1}(x),\Z),
\end{equation*}
since the $U_x$ are Stein and contractible (in a compatible way with the stratification by the topological cone theorem).
\sk
Assuming the $2$-semi-rationality (see \eqref{1}), we can prove moreover the isomorphisms
\begin{equation*}
L^2=R^2\pi'_*\Q_E/\bl(\msum_i\,\Q[E_i]_{\pi'}\br),\q L^1=R^1\pi'_*\Q_E.
\end{equation*}
These were used for the proof of the last assertion of Theorem~\ref{T1} in earlier versions of this paper.
\end{rem}

\begin{rem} \label{R2.10}
The converse of the last assertion of Theorem~\ref{T1} is not necessarily true. Indeed, there are many examples where $h^2(X)\eq h^{2n-2}(X)$, $\si(X)\eq0$, and the singularities of $X$ are locally hypersurface, isolated, and 1-semi-rational but not 2-semi-rational with $n\eq3$ using Proposition~\ref{P1}, see also Remark~\ref{R1.3}. Consider for instance the case where $X$ is locally a hypersurface of a smooth variety having only isolated singularities such that the eigenvalues of the Milnor monodromy do not contain 1 and moreover some singularity is not rational. The condition on the monodromy is equivalent to that $X$ is a $\Q$-homology manifold by \cite{Mi} (using the Wang sequence), and implies that $h^2(X)\eq h^{2n-2}(X)$ by the self-duality of intersection cohomology, see also \cite[Theorem 1]{JS}. More concretely $X$ is a hypersurface of $\PP^4$ defined by $f\eq\msum_{i=1}^4\,x_i^a\pl\msum_{i=1}^3\,x_0^{a-b}x_i^b$, where $a,b$ are mutually prime integers with $a\sgt b\gess 4$ and $a\mi b$ is equal to 1 or 5 mod 6. (The last condition is used to eliminate unwanted singularities, but it is actually sufficient to assume that $a$ is odd in the case $a{-}b$ is divisible by 3, calculating the local monodromies.) Here $X$ has a non-rational isolated singularity by \cite{exp} and \cite{new}, since $\tfrac{1}{a}\pl\tfrac{3}{b}\slt\tfrac{4}{b}\less1$. 
\sk
Even in the case where $X$ is locally a hypersurface having only {\it ordinary\one} (that is, semi-homogeneous) isolated singularities with $n\eq3$, there are many examples where the equality \eqref{3} holds and the singularities of $X$ are 1-semi-rational but not 2-semi-rational; for instance, consider the case $X\sst\PP^4$ is defined by $f\eq\msum_{i=1}^4\,x_i^a\pl\msum_{i=1}^4\,x_0^{a-b}x_i^b$ with $a\eq6$, $b\eq4$, where $h^4(X)\eq h^2(X)$ (hence $\si(X)\eq0$) according to a computation using a code in \cite{JS}. It is known that for a projective hypersurface compactification $X$ of an isolated hypersurface singularity having only one singularity, we have ${\rm def}(X)\defs h^{n+1}(X)\mi h^{n-1}(X)\eq0$ (hence $\si(X)\eq0$ when $n\eq3$) if the degree of the hypersurface $X$ is sufficiently large. This follows from \cite[Theorem 1]{JS} together with \cite[Lemma 5.3]{SS}. Note that there are also many examples with $\si(X)\eq0$ but $h^2(X)\slt h^{2n-2}(X)$ for $n\eq3$, see Remark~\ref{R2.5}.
\end{rem}

\begin{rem} \label{R2.11}
An assumption on the codimension of singular locus in \cite[Corollary 3]{JKSY} is unnecessary by using \cite[Theorem A]{KaWa} together with the Shephard-Todd-Chevalley theorem. Indeed, in the notation of that paper, we may assume that $G$ does not contain a pseudoreflection by dividing the finite group $G$ by the normal subgroup generated by pseudoreflections and taking a linear action model (and repeating this if necessary). Then there is $h\ins G$ with two eigenvalues $\al_1,\al_2$ which are different from 1 and have the {\it same order,} where the other eigenvalues of $h$ are 1 (and $\al_1,\al_2$ may coincide). Moreover the subgroup $\{g\ins G\mid L_g\eq L_h\}$ is a cyclic group, where $L_g$ denotes the fixed part by the action of $g\ins G$, and the image of $L_h$ is a local irreducible component of the {\it singular\one} locus of the quotient variety. (This does not work without assuming the nonexistence of pseudoreflections.)
\sk
It is also informed from B.~Dirks that the hypothesis on the codimension can be removed by applying \cite[Corollary 1.12]{KeSc} and \cite[Corollary 3.1]{MiVa}, since quotient singularities cannot be 1-du Bois.
\end{rem}

\section{Proof of Theorem~\ref{T2}} \label{S3}

Since $X$ is a complete intersection, it is known that $\pp\Hc^k\bl(\DD(\Z_X[n])\br)\eq0$ for $k\sgt0$, where $\pp\Hc^{\ssb}$ is the cohomology functor in \cite{BBD} and $\DD$ denotes the dual functor. (This can be shown by using the stability of semi-perversities under nearby and vanishing cycle functors, see Proposition~\ref{P4.2} below.) It is also known that Artin's vanishing theorem holds with integral coefficients, see for instance \cite{Di}, \cite[Theorem 6.0.4]{Sc}. (It was proved first in the finite and $\Z_l$-coefficient case.)
For a complete intersection $X'\sst\PP^N$ containing $X$ as a divisor (that is, $\dim X'\eq n{+}1$), we then get
\begin{equation} \label{3.1}
H^k\bl(U',\DD(\Z_{U'}[n{+}1])\br)=0\q\q\h{for any}\,\,\,\,k>0,
\end{equation}
setting $U'\defs X'\stm X$. By duality this implies
\begin{equation} \label{3.2}
H^{k+n+1}_c(U',\Z)=0\q\q\h{for any}\,\,\,\,k<0.
\end{equation}
(Note that the torsion part of $H^0\bl(U',\DD(\Z_{U'}[n{+}1])\br)$ contributes to that of $H^{k+n+2}_c(U',\Z)$.)
\sk
Using the long exact sequence associated with cohomology with compact supports, we can verify by induction on ${\rm codim}_{\PP^N}X$ that
\begin{equation} \label{3.3}
i_X^*:H^k(\PP^N,\Z)\simto H^k(X,\Z)\q\q\h{for any}\,\,\,k\les2,
\end{equation}
with $i_X\col X\into\PP^N$ the inclusion; in particular $h^1(X)=0$, $h^2(X)=1$.
The last equality implies the equivalence between condition~(b) for $X$ and that for $Y$ in the case $X$ is a projective cone of $Y$, using Remark~\ref{R2.3}.
\sk
From the commutative diagram \eqref{1.1} and the isomorphism \eqref{3.3} for $k\eq1,2$, we can deduce the vanishing
\begin{equation} \label{3.4}
H^k(X,\OO_X)=0\q\q\h{for}\,\,\,k=1,2.
\end{equation}
Considering the restriction to $X$ of the exponential sequence on $\PP^N$ and using the morphism between the associated long exact sequences, we then get the commutative diagram
\begin{equation} \label{3.5}
\begin{gathered}
\xymatrix{\Pic(\PP^N) \ar[r]^{\!\!\!\sim} \ar[d] & H^2(\PP^N,\Z) \ar[d]^{i_X^*} \\ \Pic(X) \ar[r]^{\!\!\!\!\sim} & H^2(X,\Z) }
\end{gathered}
\end{equation}
So the implication (a)$\,\Rightarrow\,$(c) follows, since condition (a) implies that ${\rm Cl}(X)\eq\Pic(X)$.
\sk
The equivalence between (c) and (d) is valid for any projective variety $X$ with $H$ fixed. Indeed, if condition~(c) holds, every irreducible reduced divisor $Z$ is linearly equivalent to a multiple of $H_X$ by a rational function $f$ on $X$ such that ${\rm div}(f)\eq Z\mi mH_X$ with $m\ins\Z_{>0}$, and this gives a defining function of $Z$ in the affine ring $A_H$. Conversely, if condition~(d) holds, every irreducible element of $A_H$ defines an irreducible divisor of $X\stm H_X$ (using Hilbert's Nullstellensatz if necessary), and every irreducible divisor of $X\stm H_X$ is contained in the zero-locus of an irreducible element $f\ins A_H$, and coincides with it using the above irreducibility of $f^{-1}(0)$, hence condition~(c) follows by viewing $f$ as a rational function on $X$. (One can see that the class of $H_X$ is not a torsion, since its $n$-fold self-intersection number is nonzero.)
\sk
Condition~(c) is clearly independent of $H$, hence so is (d) by the above argument. This gives the implication (d)$\,\Rightarrow\,$(a) using the above argument about the defining function of $Z$ in the factorial case.
Finally the equivalence between conditions (a) and (b) follows using Theorems~\ref{T1} and \ref{T4.1} below.
This completes the proof of Theorem~\ref{T2}.

\begin{rem} \label{R3.1}
It is rather easy to show that condition~(a) for $X$ implies condition~(c) for $Y$ in the case $X$ is a projective cone of $Y$ using \eqref{3.3} for $Y$ as follows: (This is related to \cite[Exer.\ II.6.3]{Ha}.) There are isomorphisms of divisor class groups
\begin{equation} \label{3.6}
{\rm Cl}(Y)\simto{\rm Cl}(X\stm\{0\})={\rm Cl}(X),
\end{equation}
since the rational function fields of $X\stm\{0\}$ and $X$ are the same. Here the first isomorphism is induced by $\rho^*$ with $\rho\col X\stm\{0\}\tos Y$ the projection of a line bundle, see for instance \cite[Proposition 1.9]{Fu}.
For an irreducible reduced divisor $Z$ on $Y$, let $Z'\sst X$ be the projective cone of $Z$. By condition~(a) for $X$, there is $f\ins\OO_{X,0}$ defining $Z'$ on a neighborhood of the vertex $0\ins X$. We have a natural $\C^*$-action $\ga_{\la}\col X\tos X$ for $\la\ins\C^*$, and $Z'\sst X$ is stable by it. This implies that $f$ can be chosen so that it is a restriction of a homogeneous polynomial of variables $x_1,\dots,x_N$, where the $x_i$ are natural homogeneous coordinates of $\PP^{N-1}$. Indeed, using complex analytic functions, we have the expansion $f\eq\msum_{m\les k<+\infty}\,f_k$ with $f_k$ restrictions to $X$ of homogeneous polynomials of degree $k$ so that $\ga_{\la}^*f_k\eq\la^kf_k$ and $f_m\nes0$. We can consider the limit of
\begin{equation*}
\la^{-m}\ga_{\la}^*f\eq\msum_{m\les k<+\infty}\,\la^{k-m}f_k\q\q\h{for}\,\,\,\,\la\to 0,
\end{equation*}
which is equal to $f_m$. Condition~(c) for $Y$ then holds using \eqref{3.6} and viewing $f_m$ as a rational function on $X$.
\end{rem}

\section{Factoriality and {\bf Q}-factoriality} \label{S4}

We can show for instance Theorem~\ref{T4.1} below, where no rationality condition is supposed, and the first assumption is satisfied in the local complete intersection case by Proposition~\ref{P4.2} below. Note also that in the latter case Theorem~\ref{T4.1} follows rather easily from \cite[X.3.4]{Gro}, see Remark~\ref{R4.4} below.

\begin{thm} \label{T4.1}
Let $X$ be a normal projective variety. Assume $\DD(\Z_X[n])\ins\pp D^{\les0}(X,\Z)$ with $\DD$ the dual functor in $D^b_c(X,\Z)$ $($see \eqref{4.9} below$)$, and the singular locus has at least codimension three. Then any effective divisor $D\ins\Div(X)$ which is $\Q$-Cartier is Cartier $($that is, Zariski-locally principal$\one)$, in particular, factoriality of $X$ is equivalent to $\Q$-factoriality.
\end{thm}

\begin{proof}
Recall that a divisor $D$ is said to be $\Q$-Cartier if $mD$ is Cartier $($that is, Zariski-locally principal$\one)$ for some positive integer $m$. 
Taking a resolution of singularities whose exceptional divisors are smooth and have normal crossings and using the direct image of an appropriate invertible sheaf together with GAGA, Nakayama's lemma, and the Hartogs-type lemma for normal schemes, the assertion is essentially analytic-local on $X$. (A similar argument is implicitly used in Remark~\ref{R2.9}.)
\sk
Let $D\eq\msum_i\,m_iD_i\ins\Div(X)$ with $D_i$ reduced, irreducible and $m_i\ins\Z_{>0}$ be an effective divisor on $X$. Let $g$ be a function on a Zariski-open subset $U\sst X$ such that ${\rm div}\,g\eq mD$ for a positive integer $m$. This gives the ramified cyclic covering
\begin{equation*}
V\defs\{t^m\eq g\}\sst U{\times}\C\to U,
\end{equation*}
with $t$ the coordinate of $\C$. Let $\Vt$ be the normalization of $V$. Then it is enough to show that $\Vt$ is \'etale over $U$ by considering the pullback of $t\,({=}\,g^{1/m})$. (In the local complete intersection case, the assertion then follows from \cite[X.3.4]{Gro}, see Remark~\ref{R4.4} below).
\sk
Let $\Sc$ be a Whitney stratification of $X$. Arguing by induction on strata and using an iterated hyperplane section transversal to a stratum, we may assume that $\Vt$ is \'etale on the complement of some point $x_0\ins X$ which is a stratum of the stratification $\Sc$. The assertion is then reduced to the vanishing
\begin{equation} \label{4.1}
H^1(X\,{\cap}\,\dd B_{x_0},\one\mu_m)=0,
\end{equation}
where $\mu_m\defs\{\la\ins\C^*\mid\la^m\eq1\}\cong\Z/m\one\Z$ and $B_{x_0}$ is a sufficiently small ball in a smooth ambient space with center $x_0$.
\sk
Since $\pp D^{\les j}$, $\pp D^{\ges j}$ are the duals of $\ppp D^{\ges-j}$, $\ppp D^{\les-j}$ (see Remark~\ref{R4.5} below), we have
\begin{equation} \label{4.2}
\Z_X[n]\in\ppp D^{\ges0}(X,\Z),
\end{equation}
that is, $\ppp\,\Hc^k(\Z_X[n])\eq0$ for any $k\slt0$.
\sk
Let $\pppo\IC_X\Z$ be the intersection complex for the dual middle perversity $\pb^+$. (It is different from the usual one even in the case of surface singularities of type $A,D,E$. Indeed, some torsion part appears from the second link cohomology and the determinant of the intersection matrix in \cite{Mu} is not $\pm1$.)
This complex can be described by iterating open direct images and truncations (see \cite{GoMP}, \cite[1.4.23 and 2.1.11]{BBD}), and we have
\begin{equation} \label{4.3}
\pppo\IC_X\Z={}^{+}\!\tau_{<0}\RR j'_*(\pppo\IC_{X'}\Z),
\end{equation}
with $X'\defs X\stm\{x_0\}$ and $j'\col X'\into X$ the inclusion.
Here ${}^{+}\!\tau_{<0}$ is the classical truncation associated with the {\it dual\one} $t$-structure on the bounded complexes of finite $\Z$-modules, see \cite[3.3.2]{BBD}. In this section we put
\begin{equation} \label{4.4}
\CCb:=C\bl(\Z_X[n]\to\pppo\IC_X\Z\br)\,\in\,\ppp D^{\ges-1}(X,\Z),
\end{equation}
(which is different from the one in Sections~\ref{S1} and \ref{S2}), where the last inclusion follows from \eqref{4.2}. Since $\dim{\rm supp}\,\CCb\less\dim\Si\less n\mi3$, we then get
\begin{equation} \label{4.5}
\Hc^{-n+k}\one\CCb=0\,\,\,\,\h{for}\,\,\,k\less1,\q\Hc^{-n+2}\one\CCb\,\,\,\h{is torsion-free,}
\end{equation}
using the truncation $^+\tau_{\les -n+1}$ as in \cite[3.3.2]{BBD}. Indeed, the above estimate of $\dim{\rm supp}\,\CCb$ implies that
\begin{equation*}
^+\tau_{\les -n+1}\CCb\ins\ppp D^{\les-2}(X,\Z),
\end{equation*}
and the canonical morphism 
\begin{equation*}
^+\Hc^{-n+1}\CCb\simto{}^+\Hc^{-n+1}\CCb
\end{equation*}
induced by $\iota:{}^+\tau_{\les -n+1}\CCb\into\CCb$ must vanish, since $\iota\eq0$ by vanishing of negative extensions:
\begin{equation*}
{\rm Hom}\bl(\ppp D^{\les-2}(X,\Z),\ppp D^{\ges-1}(X,\Z)\br)=0.
\end{equation*}
\sk
We can verify moreover that the long exact sequence associated with \eqref{4.4} implies the isomorphisms
\begin{equation} \label{4.6}
\Hc^{-n+k}(\pppo\IC_X\Z)\eq\Hc^{-n+k}\one\CCb\q\h{for any}\,\,\,k\gess1.
\end{equation}
\sk
There is the Leray-type spectral sequence
\begin{equation} \label{4.7}
E^{p,q}_2=R^p\!j'_*\Hc^q(\pppo\IC_{X'}\Z)\Longrightarrow R^{p+q}j'_*(\pppo\IC_{X'}\Z),
\end{equation}
which is associated with the usual classical truncations $\tau_{\les\ssb}$ and the usual sheaf cohomology functor $\Hc^{\ssb}$, where $R^p\!j'_*\defs\Hc^p\RR j'_*$. (Here some renumbering of spectral sequence is used, see \cite[(1.4.8)]{De}.)
\sk
In view of \eqref{4.3}, \eqref{4.5}, \eqref{4.6}, \eqref{4.7}, we conclude that
\begin{equation} \label{4.8}
E^{1,-n}_2=0,\q E^{2,-n}_2\,\,\,\h{is torsion-free.}
\end{equation}
Indeed, $E^{p,q}_2\eq0$ for any $p\ins\Z$, $q\slt{-}n$ or $q\eq{-}n{+}1$ by \eqref{4.5}, \eqref{4.6}, hence $E^{p,-n}_2\eq E^{p,-n}_r$ for any $r\gess3$ and $p\eq1,2$. Here $\Hc^{-n}(\pppo\IC_{X'}\Z)\eq\Z_{X'}[n]$ by \eqref{4.5} (and $n\gess3$ by assumption).
So the assertion follows using a universal coefficient theorem.
\end{proof}

\begin{prop} \label{P4.2}
We have the inclusion \eqref{4.2} if $X$ is a local complete intersection.
\end{prop}

\begin{proof}
It is known that the middle perversity $\bf p$ and the dual middle perversity $\pb^+$ in the $\Z$-coefficient case (see \cite[3.3]{BBD}) are stable by shifted nearby and vanishing cycle functors $\psi_h[-1],\varphi_h[-1]$ for a locally defined function $h$ on $X$, see for instance \cite[Theorem 6.0.2 and Example 6.0.2]{Sc} (and also \cite[1.1]{int} and the references noted there). By assumption $X$ is locally an intersection of hypersurfaces $h_j^{-1}(0)$ for $j\eq1,\dots,r$ in a smooth variety $Y$. Using the distinguished triangles
\begin{equation*}
i_{h_j}^*\tos\psi_{h_j}\tos\varphi_{h_j}\buildrel{\!\![1]}\over\to
\end{equation*}
with $i_{h_j}$ the inclusion of $h_j^{-1}(0)$ in $X$, we see that the shifted constant sheaf $\Z_X[n]$ is locally quasi-isomorphic to the complex obtained by iterating the shifted mapping cone
\begin{equation*}
C\bl(\psi_{h_j}[-1]\tos\varphi_{h_j}[-1]\br)[-1]\q\h{for}\,\,\,j\eq1,\dots,r
\end{equation*}
to $\Z_Y[\dim Y]$. So the assertion follows.
\end{proof}

\begin{rem} \label{R4.3}
By Grothendieck (see \cite[XI, Theorem 3.13 (ii)]{Gro}) factoriality holds in the case $X$ is a local complete intersection and ${\rm codim}_X\Sing\,X\gess4$. (Indeed, parafactoriality of $\OO_{X,x}$ implies that any invertible sheaf on ${\rm Spec}\,\OO_{X,x}\stm\{x\}$ is extendable to an invertible sheaf on ${\rm Spec}\,\OO_{X,x}$, where $x$ is not necessarily a closed point of $X$. For the proof of the coincidence of $\Div(X)$ and $\CDiv(X)$ in the case all local rings (which are inductive limits of affine rings) are factorial, it may be more intuitive to employ a (trivial) extension of Hilbert's Nullstellensatz to reduced affine varieties in order to show that an irreducible element $f$ of a local ring defines locally a reduced irreducible divisor, since in the case $(gh)^m\ins(f)$ for some $m\sgt1$, we get $g\ins(f)$ or $h\ins(f)$ using factoriality.)
\end{rem}

\begin{rem} \label{R4.4}
In the local complete intersection case, Theorem~\ref{T4.1} follows from \cite[X.3.4]{Gro} by using it instead of \eqref{4.1} in the proof of Theorem~\ref{T4.1}.
\end{rem}

\begin{rem} \label{R4.5}
Let $X$ be a reduced complex analytic space. We have the semi-perversity conditions for $\Fb\ins D^b_c(X,\Z)$ as follows:
\begin{equation} \label{4.9}
\begin{aligned}
\Fb\ins\pp D^{\les j}&\iff\Hc^ki_S^*\Fb=0\q\h{if}\,\,\,k>j\mi d_S,\\
\Fb\ins\pp D^{\ges j}&\iff\Hc^ki_S^!\Fb=0\q\h{if}\,\,\,k<j\mi d_S,
\end{aligned}
\end{equation}
where $S$ runs over any strata of a Whitney stratification with $i_S$ the inclusion of $S$ into $X$ and $d_S=\dim S$, see \cite{BBD}.
We say that $\Fb$ satisfies the perversity condition if these two semi-perversity conditions hold for $j=0$.
Note that these are not dual to each other, and there is the following dual semi-perversity conditions:
\begin{equation} \label{4.10}
\begin{aligned}
\Fb\ins\ppp D^{\les j}&\iff\Hc^ki_S^*\Fb=\begin{cases}0&\h{if}\,\,\,k>j\mi d_S\pl 1,\\ \h{torsion}&\h{if}\,\,\,k=j\mi d_S\pl 1,\end{cases}\\
\Fb\ins\ppp D^{\ges j}&\iff\Hc^ki_S^!\Fb=\begin{cases}0&\h{if}\,\,\,k<j\mi d_S,\\\h{torsion-free}&\h{if}\,\,\,k=j\mi d_S.\end{cases}
\end{aligned}
\end{equation}
These are the dual of $\pp D^{\ges -j}$ and $\pp D^{\les -j}$ respectively, see \cite[3.3]{BBD}.
The stability of these semi-perversity conditions by shifted nearby and vanishing cycle functors is shown in \cite{Sc}. (Note that the upper semi-perversity is related to rectified homotopical (or homological) depth, see for instance \cite{Gro}, \cite{HaLe}, \cite{RSW2}.)
\end{rem}

\begin{rem} \label{R4.6}
Let $(X,0)$ be a germ of a normal complex analytic space of dimension 2. Set $E\defs\Hc^0(\pppo\IC_X\Z)_0$. This is a finite abelian group, and is identified with a sheaf supported at 0 (since the direct image by a closed immersion is usually omitted). We have a self-dual distinguished triangle
\begin{equation*}
\pp\,\IC_X\Z\to\pppo\IC_X\Z\to E\buildrel{\!\![1]}\over\to,
\end{equation*}
with $\DD(\pp\,\IC_X\Z)\eq\pppo\IC_X\Z$ and $\DD(E)\eq E[-1]$.
\sk
Let $M$ be the intersection matrix of the exceptional divisors of the minimal resolution of $(X,0)$. From a calculation of a weight spectral sequence, it seems to follow that
\begin{equation} \label{4.11}
|E|=|\det(M)|.
\end{equation}
\sk
Assume $(X,0)$ is a hypersurface. Then the second link cohomology is given by $T\mi\id$ with $T$ the monodromy on the vanishing cohomology using the Wang sequence and Alexander duality (or local cohomology) as in \cite{Mi}. This should imply that
\begin{equation} \label{4.12}
|E|=|\det(T\mi\id)|.
\end{equation}
\sk
In the rational double point case (that is, of type $A,D,E$), these are compatible employing the Picard-Lefschetz formula (see for instance \cite{La}) and Brieskorn's theory of simultaneous resolutions \cite{Br}. In the case of type $A_k$ (where the eigenvalues of the monodromy are $e^{2\pi\sqrt{-1}\,j/(k+1)}$ for $j\ins\{1,\dots,k\}$), this argument implies that
\begin{equation} \label{4.13}
|\det(M)|=(x^{k+1}\mi1)/(x\mi1)|_{x=1}=k\pl1,
\end{equation} 
where $M_{i,j}\eq{-}2$ if $i\eq j$, $1$ if $|i\mi j|\eq1$, and 0 otherwise as well known, see also \cite[Appendix]{int}.
\end{rem}

\begin{rem} \label{R4.7}
It seems to be known to specialists that factoriality is equivalent to $\Q$-factoriality at least in the case of normal projective varieties of dimension not less than three having only isolated hypersurface singularities. 
An argument in Ann.\ Math.\ 127, p.~116 however seems to be misstated in an incredible way, since the cyclic covering is ramified over the divisor, and one has to take its {\it normalization,} which is \'etale over the complement of the isolated singularities.
(Here one uses Milnor's Bouquet theorem together with the Wang sequence and Alexander duality or the local cohomology for the closed submanifold $X\cap\dd B_x$ in the sphere $\dd B_x$, see \cite{Mi} where the Milnor fibration is defined on $\dd B_x\stm X$ by $f/|f|$. Note that this argument cannot be extended to the case $n\eq2$ even if 1 is not an eigenvalue of the Milnor monodromy $T$, since the variation $T\mi\id$ is not invertible with $\Z$-coefficients and some torsion appears in the second link cohomology, see Remark~\ref{R4.6}.)
\end{rem}


\begin{thebibliography}{GNPP\,88}
\bibitem[BBD\,82]{BBD} Beilinson, A., Bernstein, J., Deligne, P., Faisceaux pervers, Ast\'erisque 100 (1982).
\bibitem[dBo\,81]{dBo} du Bois, Ph., Complexe de de Rham filtr\'e d'une vari\'et\'e singuli\`ere, Bull.\ Soc.\ Math.\ France 109 (1981), 41--81.
\bibitem[BrNo\,14]{BN} Brevik, J., Nollet, S., Developments in Noether-Lefschetz theory, Contemp.\ Math., 608, A.M.S. 2014, 21--50.
\bibitem[Bri\,68]{Br} Brieskorn, E., Die Aufl\"osung der rationalen Singularit\"aten holomorpher Abbildungen, Math.\ Ann.\ 178 (1968), 255--270.
\bibitem[Che\,10]{Che} Cheltsov, I., Factorial threefold hypersurfaces, J. Alg.\ Geom.\ 19 (2010), 781--791.
\bibitem[DGPS\,23]{Sing} Decker, W., Greuel, G.-M., Pfister, G., Sch\"onemann, H., {\sc Singular} 4.3.2 --- A computer algebra system for polynomial computations, available at http://www.singular.uni-kl.de (2023).
\bibitem[De\,71]{De} Deligne, P., Th\'eorie de Hodge II, Publ.\ Math.\ I.H.E.S. 40 (1971), 5--57.
\bibitem[Di\,04]{Di} Dimca, A., Sheaves in Topology, Springer, Berlin, 2004.
\bibitem[Fu\,84]{Fu} Fulton, W., Intersection Theory, Springer, New York, 1984.
\bibitem[GoMP\,83]{GoMP} Goresky, M., MacPherson, R., Intersection homology II, Inv. Math. 72 (1983), 77--129. 
\bibitem[Gra\,62]{Gra} Grauert, H., \"Uber Modifikationen und exzeptionelle analytische Mengen, Math.\ Ann.\ 146 (1962), 331--368.
\bibitem[GrRi\,70]{GrRi} Grauert, H., Riemenschneider, O., Verschwindungss\"atze f\"ur analytische Kohomologiegruppen auf komplexen R\"aumen, Inv.\ Math.\ 11 (1970), 263--292.
\bibitem[Gro\,68]{Gro} Grothendieck, A., Cohomologie locale des faisceaux coh\'erents et th\'eor\`emes de Lefschetz locaux et globaux (SGA 2), North Holland, Amsterdam, 1968.
\bibitem[GNPP\,88]{GNPP} Guill\'en, F., Navarro Aznar, V., Pascual Gainza, P., Puerta, F., Hyperr\'esolutions cubiques et descente cohomologique, Lect.\ Notes in Math.\  1335, Springer, Berlin, 1988.
\bibitem[HaL\^e\,91]{HaLe} Hamm, H.A., L\^e, D.T., Rectified homotopical depth and Grothendieck conjectures, Grothendieck Festschrift II, Birkh\"auser, Boston, 1991, 311--351.
\bibitem[Har\,77]{Ha} Hartshorne, R., Algebraic Geometry, Springer, New York, 1977.
\bibitem[Ho\,75]{Ho} Hornell, J., A note on the geometric criteria for the factoriality of an affine ring, Proc.\ Amer.\ Math.\ Soc.\ 53 (1975), 45--50.
\bibitem[JKSY\,22]{JKSY} Jung, S.-J., Kim, I.-K., Saito, M., Yoon, Y., Higher Du Bois singularities of hypersurfaces, Proc. Lond. Math. Soc. (3) 125 (2022), 543--567.
\bibitem[JS\,25]{JS} Jung, S.-J., Saito, M., Defect of projective hypersurfaces with isolated singularities (arxiv: 2512.23522).
\bibitem[KaWa\,82]{KaWa} Kac, V., Watanabe, K.-i., Finite linear groups whose ring of invariants is a complete intersection, Bull.\ Amer.\ Math.\ Soc.\ (N.S.) 6 (1982), 221--223.
\bibitem[KeSc\,21]{KeSc} Kebekus, S., Schnell, C., Extending holomorphic forms from the regular locus of a complex space to a resolution of singularities, J. Am. Math. Soc. 34 (2021), 315--368.
\bibitem[Ko\,99]{Ko} Kov\'acs, S.J., Rational, log canonical, Du Bois singularities: on the conjectures of Koll\'ar and Steenbrink, Comp.\ Math.\ 118 (1999), 123--133.
\bibitem[La\,81]{La} Lamotke, K., The topology of complex projective varieties after S.\ Lefschetz, Topology 20 (1981), 15--51.
\bibitem[MiVa\,19]{MiVa} Miller, C., Vassiliadou, S., (Co)torsion of exterior powers of differentials over complete intersections, J. Sing. 19 (2019), 131--162.
\bibitem[Mi\,68]{Mi} Milnor, J., Singular Points of Complex Hypersurfaces, Princeton Univ.\ Press, 1968.
\bibitem[Mu\,61]{Mu} Mumford, D., The topology of normal singularities of an algebraic surface and a criterion for simplicity, Publ.\ Math.\ I.H.E.S. 9 (1961), 5--22.
\bibitem[NaSt\,95]{NaSt} Namikawa, Y, Steenbrink, J.H.M., Global smoothing of Calabi-Yau threefolds, Inv.\ Math.\ 122 (1995), 403--419.
\bibitem[PaPo\,25]{PaPo} Park, S.G., Popa, M., $\Q$-factoriality and Hodge-Du Bois theory (arxiv:2508.17748).
\bibitem[RSW\,21]{RSW} Reichelt, T., Saito, M., Walther, U., Dependence of Lyubeznik numbers of cones of projective schemes on projective embeddings, Selecta Math.\ (N.S.) 27 (2021), Paper No.\ 6, 22 pp.
\bibitem[RSW\,23]{RSW2} Reichelt, T., Saito, M., Walther, U., Topological calculation of local cohomological dimension, J. Sing. 26 (2023), 13--22.
\bibitem[Sa\,83]{exp} Saito, M., On the exponents and the geometric genus of an isolated hypersurface singularity, Proc. Symp. Pure Math. 40, Part 2, A.M.S. (1983), 465--472.
\bibitem[Sa\,88]{new} Saito, M.,Exponents and Newton polyhedra of isolated hypersurface singularities, Math.\ Ann.\ 281 (1988), 411--417.
\bibitem[Sa\,90]{mhm} Saito, M., Mixed Hodge modules, Publ. RIMS, Kyoto Univ.\ 26 (1990), 221--333.
\bibitem[Sa\,91]{kol} Saito, M., On Koll\'ar's conjecture, Proc.\ Symp.\ Pure Math., 52, Part 2, A.M.S. Providence, RI, 1991, pp~509--517.
\bibitem[Sa\,00]{mhc} Saito, M., Mixed Hodge complexes on algebraic varieties, Math.\ Ann.\ 316 (2000), 283–331.
\bibitem[Sa\,02]{bconj} Saito, M., Bloch's conjecture, Deligne cohomology and higher Chow groups (arxiv:math /9910113v14, 2002). 
\bibitem[Sa\,20]{int} Saito, M., Lowest non-zero vanishing cohomology of holomorphic functions (arxiv:2008.10529).
\bibitem[Sa\,25]{wh} Saito, M., Bernstein-Sato polynomials for projective hypersurfaces with weighted homogeneous isolated singularities (arxiv:1609.04801v11).
\bibitem[ScSt\,85]{SS} Scherk, J., Steenbrink, J.H.M., On the mixed Hodge structure on the cohomology of the Milnor fibre, Math.\ Ann.\ 271 (1985), 641--665.
\bibitem[Sc\,01]{Sc} Sch\"urmann, J., Topology of Singular Spaces and Constructible Sheaves, Monografie Matematyczne, vol.~63, Birkh\"auser, Basel, 2001.
\end{thebibliography}
\end{document}